\documentclass[12pt]{amsart}

\usepackage{times,amsfonts,amsmath,amstext,amsbsy,amssymb,
amsopn,amsthm,upref,eucal, mathrsfs}

\usepackage{url}

\def\thispartnum{III}
\input extref-arxiv.tex

\newcommand \rrad {\mathfrak{r}}

\newcommand \la {\lambda}

\DeclareMathOperator{\Span}{span}

\DeclareMathOperator{\Int}{int}

\DeclareMathOperator{\tr}{tr}

\newcommand \Mat {{\mathrm{Mat}}}
\newcommand \Matreg{{\mathrm{Mat}_{\mathrm{reg}}}}

\newcommand \Prob {{\mathbb P}}

\DeclareMathOperator{\Leb}{Leb}

\newcommand  \const {\mathrm{const}}
\newcommand  \loc {{\mathrm{loc}}}

\newcommand  \fin {{\mathrm{fin}}}
\newcommand  \Mfin {\mathfrak{M}_\fin}

\newcommand \Conf {{\mathrm {Conf}}}

\newcommand \conf {{\mathrm {conf}}}

\newcommand\scrI{{\mathscr I}}

\newtheorem{theorem}{Theorem}[section]

\newtheorem{lemma}[theorem]{Lemma}

\newtheorem{corollary}[theorem]{Corollary}
\newtheorem{proposition}[theorem]{Proposition}

\begin{document}
\title[The Ergodic Decomposition of Infinite Pickrell Measures. III]{Infinite Determinantal Measures and
The Ergodic Decomposition of Infinite Pickrell Measures III. The infinite Bessel process as the limit of radial parts of finite-dimensional projections of infinite Pickrell measures}

\author{Alexander I. Bufetov}

\address{ Aix-Marseille Universit{\'e}, Centrale Marseille, CNRS,   I2M}
\address{UMR 7373, rue F. Joliot Curie , Marseille, France}
\address{Steklov Institute of Mathematics,
Moscow, Russia}
\address{Institute for Information Transmission Problems,
 Moscow, Russia}
\address{National Research University Higher School of Economics,
 Moscow, }
 \address{Russia}

\begin{abstract} The third part of the paper concludes the proof of the main result --- the 
description of the ergodic decomposition of infinite Pickrell measures.
 First it is shown that the scaling limit of radial parts of finite-dimensional 
 infinite Pickrell measures is precisely the infinite Bessel point process. It is then established that 
 the  ``gaussian parameter'' almost surely vanishes for our ergodic components, and the convergence to the scaling limit is then established in the space of finite measures on the space of finite measures. Finally, singularity is established for Pickrell measures corresponding to different values of the parameter. 
\end{abstract}

\maketitle
\tableofcontents
\section{Introduction}

This paper is the third and final of the cycle of three papers giving the explicit construction of the ergodic decomposition of infinite Pickrell measures.
Quotes to the other parts of the paper \cite{infdet1, infdet2} are organized as follows: Proposition II.2.3, equation (I.9), etc.

The paper is organized as follows.
In Section 2, we go back to radial parts of Pickrell measures.
We start by recalling the  determinantal representation for
radial parts of finite Pickrell measures and the convergence of the resulting determinantal processes to the modified Bessel point process (the usual  Bessel point process of Tracy and Widom \cite{TracyWidom} subject to the change of variable $y=4/x$). Next, we represent the radial parts of infinite Pickrell measures as infinite determinantal measures corresponding to finite-rank perturbations
of Jacobi orthogonal polynomial ensembles. The main result of this section is Proposition \ref{conv-bns-bs} which shows that the scaling limit of the infinite determinantal measures corresponding to the radial parts of infinite Pickrell measures  is precisely the modified infinite Bessel point process of the Introduction.
Infinite determinantal measures are made finite by taking the product with a suitable multiplicative functional, and weak convergence is established both in the space of finite measures on the space of configurations and in the space of finite measures in the space of finite measures. The latter statement will be essential in the proof of the vanishing of the ``Gaussian parameter'' in the following section.

In Section 3, we pass from the convergence, in the space of finite measures on the space of configurations and in the space of finite measures in the space of finite measures,  of rescaled radial parts of Pickrell
measures to the convergence, in the space of finite measures on the Pickrell set, of finite-dimensional approximations of Pickrell measures. In particular, in this section we establish the vanishing of the ``Gaussian parameter'' for ergodic components of infinite  Pickrell measures. Proposition \ref{weakconfpick} proved in this section allows us  to  complete the proof of Proposition \ref{weak-pick}.

The final Section 4 is mainly devoted to the proof of Lemma \ref{main-lemma}, which relies on the well-known asymptotics of the Harish-Chandra-Itzykson-Zuber orbital integrals. Combining Lemma \ref{main-lemma} with Proposition \ref{weak-pick}, we conclude the proof of Theorem \ref{mainthm}. Using Kakutani's Theorem in the same way as Borodin and Olshanski \cite{BO}, we prove that Pickrell measures corresponding to distrinct values of the parameter $s$ are mutually singular.

\section{Weak convergence of rescaled radial parts of Pickrell measures}

\subsection{ The case $s>-1$: finite Pickrell measures}
\subsubsection{Determinantal representation of the radial parts of  finite Pickrell measures}

We go back to radial parts of Pickrell measures $\mu_n^{(s)}$ and start with the case $s>-1$ .
Recall that $P_n^{(s)}$ stand for the Jacobi polynomials corresponding to the weight $(1-u)^s$ on the interval $[-1,1]$.

We start by giving a determinantal representation for the radial part of finite Pickrell measures: in other words, we simply rewrite the formula (\ref{det-CDJ}) in the
coordinates $\la_1, \dots, \la_n$.
Set
\begin{multline*}
{\hat K}_n^{(s)} (\la_1, \la_2) = \frac{n(n+s)}{2n+s} \frac{1}{(1+\la_1)^{s/2}
(1+\la_2)^{s/2}}\times \\
\times \frac{
P_n^{(s)} \left( \frac{\la_1-1}{\la_1+1} \right)
P_{n-1}^{(s)} \left( \frac{\la_2-1}{\la_2+1} \right) -
P_n^{(s)} \left( \frac{\la_2-1}{\la_2+1} \right)
P_{n-1}^{(s)} \left( \frac{\la_1-1}{\la_1+1} \right)
}{\la_1-\la_2}.
\end{multline*}
The kernel ${\hat K}_n^{(s)}$ is the image of the Christoffel-Darboux kernel  ${\tilde K}^{(s)}_n$ (cf. (\ref{CDJ})) under the change of variable $$u_i=\frac{\lambda_i-1}{\lambda_i+1}.$$
Another representation for the kernel ${\hat K}_n^{(s)}$ is
\begin{multline*}
{\hat K}_n^{(s)} (\la_1, \la_2) = \frac{1}{(1+\la_1)^{s/2+1}(1+\la_2)^{s/2+1}}\times\\
\times{\displaystyle\sum\limits_{l=0}^{n-1}(2l+s+1){P_l^{(s)}\left( \frac{\la_1-1}{\la_1+1} \right)\cdot P_l^{(s)} \left( \frac{\la_2-1}{\la_2+1} \right)}}.
\end{multline*}

The kernel ${\hat K}_n^{(s)}$ is by definition the kernel of the operator of orthogonal projection in $L_2((0, +\infty), \Leb)$ onto the subspace
\begin{multline*}
{\hat L}^{(s, n)}=\Span \left(\frac 1{(\la+1)^{s/2+1}}P_l^{(s)}\left( \frac{\la-1}{\la+1} \right), l=0, \dots, n-1 \right)=\\
=\Span \left(\frac 1{(\la+1)^{s/2+1}}\left( \frac{\la-1}{\la+1} \right)^l, l=0, \dots, n-1 \right).
\end{multline*}

Proposition \ref{rad-jac} implies the following determinantal representation of the radial part of the Pickrell measure.
\begin{proposition}
For $s>-1$,  we have
\begin{equation*}
(\mathfrak{rad}_{n})_*\mu_n^{(s)} = \frac{1}{n!}\det{\hat  K}_n^{(s)}(\la_i,\la_j)\prod\limits_{i=1}^n d\la_i.
\end{equation*}
\end{proposition}

\subsubsection{Scaling}

For $\beta>0$, let $\mathrm{hom}_{\beta}: (0, +\infty)\to (0, +\infty)$ be the
homothety map that sends $x$ to $\beta x$; we keep the same symbol for the induced scaling transformation of  $\Conf((0, +\infty))$.

We now give an explicit determinantal representation for the measure
\begin{equation*}
\left (\conf\circ \mathrm{hom}_{n^2} \circ {\mathfrak{rad}}_n\right)_*\mu_n^{(s)},
\end{equation*}
the push-forward to the space of configurations
of the
rescaled radial part
of the Pickrell measure $\mu_n^{(s)}$.

 Consider the rescaled Christoffel-Darboux kernel
\begin{equation*}
K_n^{(s)}(\la_1,\la_2)=n^2\hat K_n^{(s)}\left( n^2 \la_1, n^2 \la_2\right)
\end{equation*}
of orthogonal projection onto the rescaled subspace
\begin{multline*}
{L}^{(s, n)}=\Span \left(\frac 1{(n^2\la+1)^{s/2+1}}P_l^{(s)}\left( \frac{n^2\la-1}{n^2\la+1} \right) \right)=
\\
=\Span \left(\frac 1{(n^2\la+1)^{s/2+1}}\left( \frac{n^2\la-1}{n^2\la+1} \right)^l, l=0, \dots, n-1 \right).
\end{multline*}
The kernel $K_n^{(s)}$ induces a determinantal process $\Prob_{K_n^{(s)}}$ on the space $\Conf((0, +\infty))$.

\begin{proposition}
For $s>-1$, we have
$$
\left(\mathrm{hom}_{n^2}\circ \mathfrak{rad}_{n}\right)_*\mu_n^{(s)} = \frac{1}{n!}\det{ K}_n^{(s)}(\la_i,\la_j)\prod\limits_{i=1}^n d\la_i.
$$
Equivalently,
$$
\left (\conf\circ \mathrm{hom}_{n^2} \circ {\mathfrak{rad}}_n\right)_*\mu_n^{(s)}=\Prob_{K_n^{(s)}}.
$$
\end{proposition}

\subsubsection{Scaling limit}

The scaling limit for radial parts of finite Pickrell measures is a variant of the well-known result of Tracy and Widom \cite{TracyWidom} claiming that the scaling limit
of Jacobi orthogonal polynomial ensembles is the Bessel point process.

\begin{proposition}
For any $s>-1$, as $n\to\infty$, the kernel $K_n^{(s)}$
converges
to the modified Bessel kernel $J^{(s)}$ uniformly in the totality of variables
on compact subsets of
$(0, +\infty)\times (0, +\infty)$.
We therefore have
$$
K_n^{(s)}\to J^{(s)}\text{ in }\scrI_{1, \loc}((0, +\infty), \Leb)
$$
and
$$
\Prob_{K_n^{(s)}}\to \Prob_{J^{(s)}}\text{ in }\Mfin\Conf((0, +\infty)).
$$
\end{proposition}

Proof. This is an immediate corollary of the classical Heine-Mehler asymptotics for Jacobi polynomials, see, e.g., Szeg{\"o} \cite{Szego}.

{\bf Remark.} As the Heine-Mehler asymptotics show, the uniform
convergence in fact takes place on arbitrary simply connected compact subsets of
$({\mathbb C}\setminus 0)\times ({\mathbb C}\setminus 0)$.

\subsection{The case $s\leq -1$: infinite Pickrell measures}
\subsubsection{Representation of radial parts of infinite Pickrell measures as infinite determinantal measures}

Our first aim is to show that for $s\leq -1$, the measure (\ref{rad-u-coord}) is an infinite determinantal measure. Similarly to the definitions given in the Introduction, set
\begin{multline*}
\hat V^{(s,n)}=\Span \biggl( \frac 1{(\la+1)^{s/2+1}},
\frac 1{(\la+1)^{s/2+1}}\left( \frac{\la-1}{\la+1} \right),  \dots,  \\
   \dots,     \frac 1{(\la+1)^{s/2+1}}P_{n-n_s}^{(s+2n_s-1)}\left( \frac{\la-1}{\la+1} \right)\biggr).
\end{multline*}

\begin{equation*}
{\hat H}^{(s,n)}=\hat V^{(s,n)}\oplus {\hat L}^{(s+2n_s, n-n_s )}.
\end{equation*}

Consider now the rescaled subspaces
\begin{multline*}
 V^{(s,n)}=\Span\biggl( \frac 1{(n^2\la+1)^{s/2+1}},
\frac 1{(n^2\la+1)^{s/2+1}}\left( \frac{n^2\la-1}{n^2\la+1} \right),  \dots,\\
  \dots,      \frac 1{(n^2\la+1)^{s/2+1}}P_{n-n_s}^{(s+2n_s-1)}\left( \frac{n^2\la-1}{n^2\la+1} \right)\biggr).
\end{multline*}

\begin{equation*}
{H}^{(s,n)}= V^{(s,n)}\oplus L^{(s+2n_s, n-n_s )}.
\end{equation*}

\begin{proposition} \label{rescaled-inf-s}Let $s\leq -1$, and let $R>0$ be arbitrary.
The radial part of the Pickrell measure is then an infinite determinantal measure
corresponding to the subspace $H={\hat H}^{(s,n)}$ and the subset $E_0=(0, R)$:
$$
( \mathfrak{rad}_{n})_*\mu_n^{(s)}={\mathbb B}\bigl({\hat H}^{(s,n)}, (0, R)\bigr).
$$
For the rescaled radial part, we have
$$\conf_*\rrad^{(n)}(\mu^{(s)})=
\left (\conf\circ \mathrm{hom}_{n^2} \circ {\mathfrak{rad}}_n\right)_*\mu_n^{(s)}={\mathbb B}\left({ H}^{(s,n)}, (0, R)\right).
$$
\end{proposition}

\subsection{The modified Bessel point process as the scaling limit of the radial parts of infinite Pickrell measures: formulation of Proposition \ref{conv-bns-bs}}

Denote ${\mathbb B}^{(s,n)}= {\mathbb B}\left({ H}^{(s,n)}, (0, R)\right)$.
We now describe the limit transition of the measures ${\mathbb B}^{(s,n)}$ to ${\mathbb B}^{(s)}$: namely, we multiply our sequence of infinite measures by a convergent multiplicative functional and
establish the convergence of the resulting sequence of determinantal probability measures.
It will be convenient to take $\beta>0$ and set $g^{\beta}(x)=\exp(-\beta x)$, while for $f$
it will be convenient to take the function $f(x)=\min(x,1)$.
Set, therefore,
$$
L^{(n, s, \beta)}=\exp(-\beta x/2) H^{(s,n)}.
$$
It is clear by definition that $L^{(n, s, \beta)}$ is a closed subspace of $L_2((0, +\infty), \Leb)$; let
$\Pi^{(n, s, \beta)}$ be the corresponding orthogonal projection operator. Recall also from (\ref{expaa}), (\ref{pisb}) the operator $\Pi^{( s, \beta)}$ of orthogonal projection onto
the subspace $L^{(s, \beta)}=\exp(-\beta x/2) H^{(s)}$.

\begin{proposition} \label{conv-bns-bs}
\begin{enumerate}
\item For all $\beta>0$ we have  $\Psi_{g^{\beta}}\in L_1\left(\Conf(0,+\infty),\mathbb{B}^{(s)}\right)$ and, for all  $n>-s+1$ we also have $\Psi_{g^{\beta}}\in L_1\left(\Conf(0,+\infty),\mathbb{B}^{(s,n)}\right)$;
\item  we have
$$\frac{\Psi_{g^{\beta}}\mathbb {B}^{(s,n)}} {\displaystyle \int\Psi_{g^{\beta}}\,d\mathbb{B}^{(s,n)}}=\Prob_{{\Pi}^{(n, s, \beta)}};$$
$$\frac{\Psi_{g^{\beta}}\mathbb {B}^{(s)}} {\displaystyle \int\Psi_{g^{\beta}}\,d\mathbb {B}^{(s)}}=\Prob_{{\Pi}^{(s, \beta)}};$$
\item We have
$$\Pi^{(n, s, \beta)}\to\Pi^ {s, \beta}\;\text{ in }\; \scrI_{1,\loc}((0, +\infty), \Leb))
\;\text{ as }\; n\to\infty\;$$
and, consequently,
$$\Prob_{{\Pi}^{(n, s, \beta)}}\to\Prob_{{\Pi}^{(s, \beta)}}$$
as $n\to\infty$ weakly in $\Mfin\bigl(\Conf((0,+\infty))\bigr)$;
\item for $f(x)=\min(x,1)$ we have $$\sqrt{f}\Pi^{(n, s, \beta)}\sqrt{f}, \sqrt{f}\Pi^ {s, \beta}\sqrt{f}\in \scrI_{1}((0, +\infty), \Leb));$$
$$\sqrt{f}\Pi^{(n, s, \beta)}\sqrt{f}\to\sqrt{f}\Pi^ {s, \beta}\sqrt{f}\;\text{ in }\; \scrI_{1}((0, +\infty), \Leb))
\;\text{ as }\; n\to\infty\;$$
and, consequently,
$$(\sigma_f)_*\Prob_{{\Pi}^{(g,s,n)}}\to(\sigma_f)_*\Prob_{{\Pi}^{(g,s)}}$$
as $n\to\infty$ weakly in $\Mfin\left(\Mfin\bigl((0,+\infty)\bigr)\right)$.
\end{enumerate}
\end{proposition}
The proof  of Proposition \ref{conv-bns-bs} will occupy the remainder of this section.

\subsection{Proof of Proposition \ref{conv-bns-bs}}

\subsubsection{Proof of the first three claims}

For $s>-1$,  write
$$
L^{(n, s, \beta)}=\exp(-\beta x/2) L^{(s,n)}_{Jac}, \ L^{(s, \beta)}=\exp(-\beta x/2) L^{(s)}
$$
and keep the notation
$\Pi^{(n, s, \beta)}$ , $\Pi^{( s, \beta)}$  for
the corresponding orthogonal projection operators.
For $s>-1$, using the Proposition \ref{conv-induced} on the convergence of induced processes, we clearly have
$$\frac{\Psi_{g^{\beta}}\Prob_{K^{(s)}_n}} {\displaystyle \int\Psi_{g^{\beta}}\,d\Prob_{K^{(s)}_n}}=\Prob_{{\Pi}^{(n, s, \beta)}};$$
$$\frac{\Psi_{g^{\beta}}\Prob_{J^{(s)}}} {\displaystyle \int\Psi_{g^{\beta}}\,d\Prob_{J^{(s)}}}=\Prob_{{\Pi}^{(s, \beta)}},$$
and also
$$\Pi^{(n, s, \beta)}\to\Pi^ {s, \beta}\text{ in }\scrI_{1,\loc}((0, +\infty), \Leb))
\;\text{ as }\; n\to\infty\;.$$

If $x_n\to x$ as $n\to\infty$, then, of course, for any $\alpha\in {\mathbb R}$ we have
$$
\lim\limits_{n\to\infty}\frac 1{n^{2\alpha}} (n^2x_n+1)^{\alpha}=x^{\alpha},
$$
and, by the Heine-Mehler classical asymptotics, for any $\alpha>-1$, we also have
$$
\lim\limits_{n\to\infty} \frac{1}{(n^2x+1)^{\alpha/2+1}}P_{n}^{(\alpha)}\left( \frac{n^2x_n-1}{n^2x_n+1} \right)=\frac{J_{\alpha}(2/{\sqrt{x}})}{\sqrt{x}}.
$$

We note the following statement on linear independence, which is immediate from Proposition \ref{lin-ind} by the change of variables $y=4/x$.

\begin{proposition}For any $s\leq -1$, and any $R>0$ the functions
\begin{equation*}
x^{-s/2-1}\chi_{(0,R)},\dots,\frac{J_{s+2n_s-1}(\frac{2}{\sqrt{x}})}{\sqrt{x}}\chi_{(0,R)}
\end{equation*}
are linearly independent and, furthermore,
are independent from the space $\chi_{(0, R)}L^{s+2n_s}$.
\end{proposition}

The proof of Proposition \ref{lin-ind} also gives, of course, that the functions
\begin{equation*}
e^{-\beta x/2}x^{-s/2-1},\dots,e^{-\beta
x/2}\frac{J_{s+2n_s-1}(\frac{2}{\sqrt{x}})}{\sqrt{x}}
\end{equation*}
are linearly independent and, furthermore, independent from the space $e^{-\beta x/2}L^{s+2n_s}$.
The first three claims of Proposition \ref{conv-bns-bs} follow now from its abstract counterparts established in the previous subsections: the first and the second claim follow from Corollary \ref{fin-rank-mult}, while
the third claim, from Proposition \ref{conv-inf-det}.
We proceed to the proof of the fourth and last claim of Proposition \ref{conv-bns-bs}.

\subsubsection{The asymptotics of $J^{(s)}$  at $0$ and at $\infty$}

We shall need the asymptotics of the modified Bessel kernel $J^{(s)}$ at $0$ and at $\infty$.

We recall that the Bessel function is denoted by $J_s$, and the usual Bessel kernel is denoted by ${\tilde J}_s$. We start with a simple estimate for ${\tilde J}_s$.
\begin{proposition}  For any $s>-1$ and any
$R>0$ we have
\begin{equation}\label{Rxest}
\displaystyle \int\limits_R^{+\infty} \frac{{\tilde J}_s (y,y)}{y} dy< +\infty.
\end{equation}
\end{proposition}
Proof. Rewrite (\ref{Rxest}) in the form
$$
\int\limits_R^{+\infty} \frac1y \int\limits_0^{1} (J_s(\sqrt{ty}))^2 dt dy =
\int\limits_0^{1} dt \int\limits_{tR}^{+\infty} \frac{(J_s(\sqrt{y}))^2}{y} dy =
$$
$$
=\int\limits_0^{+\infty} \min\Bigl(\frac{y}{R}, 1\Bigr) \cdot \frac{J_s(\sqrt{y})^2}{y} dy=
\frac{1}{R}\int\limits_0^{R} J_s(\sqrt{y})^2 dy +\int\limits_R^{+\infty}
\frac{J_s(\sqrt{y})^2}{y} dy.
$$
 It is immediate from the asymptotics of the Bessel functions at zero and at
 infinity that both integrals converge, and the proposition is proved. Effectuating the charge
 of variable $y=4/x$, we arrive at the following

 \begin{proposition}
 For any $s> -1$ and any $\varepsilon>0$ there exists $\delta>0$
 such that
 $$
 \int_0^{\delta} x J^{(s)} (x,x) dx < \varepsilon.
 $$
 \end{proposition}

 We also need the following
 \begin{proposition}
 For any $R>0$ we have
 $$
 \int\limits_0^R  {\tilde J}_{s} (y,y) dy < \infty.
 $$  \end{proposition}

 Proof. First note that
 $$
 \int\limits_0^R \left(J_s(\sqrt{y})\right)^2 dy<+\infty
 $$
 since for a fixed $s>-1$ and all sufficiently small $y>0$ we have $$\left(J_s(\sqrt{y})\right)^2=O(y^s).$$
 Now, write
 $$
 \int\limits_0^R  {\tilde J}_{s} (y,y) dy= \int\limits_0^1\int\limits_0^R  {J}_{s} (\sqrt{ty})^2 dydt\leq
 $$
 $$
 \leq (R+1)\int\limits_0^R (J_s(\sqrt{y}))^2 dy<+\infty,
 $$
 and the proposition is proved. Making the change of variables $y=4/x$, we obtain

 \begin{proposition}
 For any $R>0$ we have
 $$
 \int_R^{\infty}  J^{(s)} (x,x) dx < \infty.
 $$  \end{proposition}

\subsubsection{Uniform in $n$ asymptotics at infinity  for the kernels $K^{(n, s)}$}

We  turn to the uniform asymptotic at infinity for the kernels
$K^{(n, s)}$ and the limit kernel $J^{(s)}$. This uniform asymptotic is needed to establish the last claim of Proposition \ref{conv-bns-bs}.

\begin{proposition} \label{unif-kns-infty} For any $s>-1$ and any $\varepsilon >0$ there exists
$R>0$ such that
\begin{equation}\label{knsj-inf}
\sup_{n\in \mathbb{N}} \displaystyle \int\limits_R^{+\infty} K^{(n, s)}(x,x) dx
<\varepsilon,
\end{equation}
\end{proposition}
Proof.
We start by verifying the desired estimate~(\ref{knsj-inf}) for $s>0$.
But if $s>0$ then the classical inequalities for Borel functions
and Jacobi polynomials (see e.g. Szeg{\"o} \cite{Szego}) imply the existence of a
constant $C>0$ such that for all $x\geq 1$ we have:
$$
 \sup\limits_{n\in \mathbb{N}} K^{(n, s)} (x, x)<\displaystyle \frac{C}{x^2}.
$$
The proposition for $s>0$ is now immediate.

To consider the remaining case $s\in (-1, 0]$, we recall that the kernels $K^{(n, s)}$ are rank-one perturbations of the kernels $K^{(n-1, s+2)}$ and  note the following immediate general
\begin{proposition}
Let $K_n, K, {\check K}_n, {\check K}\in \scrI_{1, \loc}((0, +\infty), \Leb)$ be locally trace-class projections acting in
$L_2((0, +\infty), \Leb)$.
Assume
\begin{enumerate}
\item $K_n\to K$, ${\check K}_n \to {\check K}$ in $\scrI_{1, \loc}((0, +\infty), \Leb)$ as $n\to\infty$;
\item for any $\varepsilon>0$ there exists $R>0$ such that
$$
\sup\limits_{n\to\infty} \tr \chi_{(R, +\infty)}K_n\chi_{(R, +\infty)}<\varepsilon, \
 \tr \chi_{(R, +\infty)}K\chi_{(R, +\infty)}<\varepsilon;
$$
\item there exists $R_0>0$ such that
$$
\tr \chi_{(R_0, +\infty)}{\check K}\chi_{(R_0, +\infty)}<\varepsilon;
$$
\item the projection operator ${\check K}_n$ is a rank one perturbation of $K_n$.
\end{enumerate}
Then  for any $\varepsilon>0$ there exists $R>0$ such that
$$
\sup\limits_{n\to\infty} \tr \chi_{(R, +\infty)}{\check K}_n\chi_{(R, +\infty)}<\varepsilon.
$$
\end{proposition}

Proposition \ref{unif-kns-infty} is now proved completely.

\subsubsection{Uniform in $n$ asymptotics at zero for the kernels $K^{(n, s)}$ and completion of the proof of Proposition \ref{conv-bns-bs}}

We next  turn to the uniform asymptotics at zero for the kernels
$K^{(n, s)}$ and the limit kernel $J^{(s)}$. Again, this uniform asymptotics is needed to establish the last claim of Proposition \ref{conv-bns-bs}.

\begin{proposition} \label{unifest-zero} For any $\varepsilon>0$ there exists $\delta>0$ such
that for all $n \in \mathbb{N}$ we have
\begin{equation*}
\int_0^{\delta} x{K}^{(n, s)}(x,x) dx < \varepsilon.
\end{equation*}
\end{proposition}

Proof.  Going back to the $u$-variable, we reformulate our proposition as follows:
\begin{proposition}
For any $\varepsilon > 0$ there exists $R > 0$, $n_0 \in \mathbb{N}$, such that for all $n > n_0$ we have
\begin{equation*}
\frac{1}{n^2} \int\limits_{-1}^{1 - {R}/{n^2}} \frac{1+u}{1-u} \tilde{K}_n^{(s)}(u,u) du < \varepsilon.
\end{equation*}
\end{proposition}
First note that the function $\frac{1 + u}{1 - u}$ is bounded above on $[-1,0]$, and therefore
$$
\frac{1}{n^2} \int\limits_{-1}^0 \frac{1+u}{1-u} \tilde{K}_n^{(s)}(u,u)du \leq \frac{2}{n^2} \int\limits_{-1}^1 \tilde{K}_{n}^{(s)}(u,u) du = \frac{2}{n}.
$$
We proceed to estimating
$$
\frac{1}{n^2} \int\limits_0^{1-{R}/{n^2}} \frac{1 + u}{1 - u}\tilde{K}_n^{(s)}(u,u) du.
$$

Fix $\kappa>0$ (the precise choice of $\kappa$ will be described later).

Write
\begin{multline*}
\tilde{K}_n^{(s)}(u,u) = \left( \sum\limits_{l\leq \kappa n} (2l + s + 1) \bigl( P_l^{(s)}(u)\bigr)^2 \right) (1-u)^s +\\+ \left( \sum\limits_{l>\kappa n}(2l+s+1) \bigl( P_l^{(s)}(u) \bigr)^2 \right) (1-u)^s
\end{multline*}
We start by estimating
\begin{equation} \label{sharp-edge}
\frac{1}{n^2} \sum\limits_{ l \leq \kappa n} \int\limits_{1 - {1}/{l^2}}^{1 - {R}/{n^2}}\frac{1 +u}{1 - u}(2l + s + 1) \bigl(P_{l}^{(s)}(u)\bigr)^2 (1-u)^s du
\end{equation}

Using the trivial estimate
$$
\max\limits_{u \in [-1,1]} \left| P_{l}^{(s)}(u) \right| = O(l^2)
$$
we arrive, for the integral (\ref{sharp-edge}), at the upper bound
\begin{equation}\label{se1}
\const \cdot \frac{1}{n^2} \sum\limits_{l \leq \kappa n} l^{2s+1} \cdot \int\limits_{1-\frac{1}{l^2}}^{1- \frac{c}{n^2}} (1 - u)^{s-1} du
\end{equation}

We now consider three cases: $s>0$, $s = 0$, and $-1 < s < 0$.

{\it The First Case.}
If $s >0$, then the integral (\ref{se1}) is estimated above by the expression
$$
\const \cdot \frac{1}{n^2} \cdot \sum\limits_{l \leq \kappa n} l^{2s + 1}\frac{1}{l^{2s}} \leq \const \cdot \kappa^2.
$$

{\it The Second Case.} If $s = 0$, then the integral (\ref{se1}) is estimated above by the expression
$$
\const \cdot \frac{1}{n^2} \sum\limits_{l \leq \kappa n} l \cdot \log(n/l) \leq \const \cdot \kappa^2.
$$

{\it The Third Case.}
Finally, if $-1 < s < 0$, then we arrive, for the integral (\ref{se1}), at the upper bound
$$
\const \cdot \frac{1}{n^2} \left( \sum\limits_{l \leq \kappa n} l^{2s+1} \right) \cdot R^s n^{-2s} \leq \const \cdot R^s \kappa^{2 + 2s}
$$

Note that in this case, the upper bound decreases as $R$ grows.
Note that in all three cases the contribution of the integral (\ref{se1}) can be made arbitrarily small by choosing $\kappa$ sufficiently small. We next estimate
\begin{equation} \label{se2}
\frac{1}{n^2} \sum\limits_{l \leq \kappa n} \int\limits_{0}^{1-\frac{1}{l^2}} \frac{1 +u}{1 - u}(2l + s + 1) \left( P_l^{(s)}(u) \right)^2 (1-u)^s du
\end{equation}
Here we use the estimate (7.32.5) in Szeg{\"o} \cite{Szego} that gives
$$
\left| P_l^{(s)}(u) \right| \leq \const \frac{(1-u)^{-\frac{s}{2}-\frac{1}{4}}}{\sqrt{l}}
$$
as long as $u \in [0, 1 -\frac{1}{l^2}]$ and arrive, for the integral (\ref{se2}), at the upper bound
$$
\const \cdot \frac{1}{n^2} \sum\limits_{l \leq \kappa n} l \leq \const \cdot \kappa^2
$$
which, again, can be made arbitrarily small as soon as $\kappa$ is chosen sufficiently small.

It remains to estimate the integral
\begin{equation} \label{se3}
\frac{1}{n^2} \sum\limits_{\kappa n \leq l < n} \int\limits_{0}^{1 - \frac{R}{n^2}} \frac{1 + u}{1 - u} \cdot (2l + s + 1) (P_l^{(s)})^2 (1-u)^s du
\end{equation}

Here again we use the estimate (7.32.5) in Szego \cite{Szego} and note that since the ratio $l/n$ is bounded below, we have a uniform estimate
$$
\left| P_l^{(s)}(u) \right| \leq \const \cdot \frac{(1-u)^{-\frac{s}{2} -\frac{1}{4} }}{\sqrt{l}}
$$
valid as long as $\kappa n \leq l \leq n$, $u \in [0, 1 -\frac{R}{n^2}]$, and in which the constant depends on $\kappa$ and does not grow as $R$ grows.

We thus arrive, for integral (\ref{se3}), at the upper bound
$$
\frac{\const}{n^3} \sum\limits_{\kappa n \leq l < n} \int\limits_{0}^{1 - \frac{R}{n^2}} (1-u)^{-\frac{3}{2}}
 du \leq \frac{\const}{n^2}\left(1-\frac{n}{\sqrt{R}}\right).$$
Now choosing $\kappa$ sufficiently small as a function of $\varepsilon$ and then $R$ sufficiently large as a function of $\varepsilon$ and $\kappa$, we conclude the proof of the proposition.

The fourth claim of Proposition \ref{conv-bns-bs} is now an immediate corollary of uniform estimates given in Propositions \ref{unif-kns-infty}, \ref{unifest-zero} and the general statement given in Proposition \ref{conv-unif-finrank}.

Proposition \ref{conv-bns-bs} is proved completely.

\section{Convergence of approximating measures on the Pickrell set and proof of Propositions \ref{expbetaint},  \ref{weak-pick}}

\subsection{Proof of Proposition \ref{expbetaint}}

Proposition \ref{expbetaint}  easily follows  from what has already been established.
Recall that we have a natural forgetting map
$\conf:\Omega_P\to \Conf(0,+\infty)$ that assigns
to $\omega=(\gamma, x)$, $x=(x_1, \dots, x_n, \dots)$, the configuration
$\omega(x)=(x_1, \dots, x_n, \dots)$.
By definition, the map $\conf$ is $\rrad^{(n)}(\mu^{(s)})$-almost surely bijective. The characterization of the measure
$\conf_*\rrad^{(n)}(\mu^{(s)})$ as an infinite determinantal measure given by
Proposition \ref{rescaled-inf-s} and  the first statement of Proposition \ref{conv-bns-bs} now imply Proposition \ref{expbetaint}.
We proceed to the  proof of  Proposition \ref{weak-pick}.

\subsection{Proof of  Proposition \ref{weak-pick}}
Recall that, by definition,  we have
$$
\conf_*\nu^{(s,n, \beta)}=\Prob_{\Pi^{(s, n, \beta)}}.
$$
Recall that Proposition \ref{conv-bns-bs} implies
that, for any $s\in {\mathbb R}$, $\beta>0$, as $n\to\infty$ we have
$$
\Prob_{\Pi^{(s, n, \beta)}}\to \Prob_{\Pi^{(s,  \beta)}}
$$
in $\Mfin(\Conf((0, +\infty)))$ and,
furthermore,  setting $f(x)=\min(x,1)$, also the weak convergence
 $$
(\sigma_f)_*\Prob_{\Pi^{(s, n, \beta)}}\to (\sigma_f)_*\Prob_{\Pi^{(s,  \beta)}}
$$
in $\Mfin(\Mfin((0, +\infty)))$.
We  now  need to pass from
weak convergence of probability measures on the space of configurations established in
Proposition \ref{conv-bns-bs} to the weak convergence of probability measures on the Pickrell set.

We have a natural map
$${\mathfrak s}\colon\Omega_P\to\Mfin\bigl((0,+\infty)\bigr)$$
defined by the formula
$${\mathfrak s}(\omega)=\sum\limits_{i=1}^{\infty}\min (x_i(\omega), 1)\delta_{x_i(\omega)}\,.$$

The map ${\mathfrak s}$ is bijective in restriction to the subset $\Omega_P^0$ defined, we recall,
as the subset of $\omega=(\gamma, x)\in\Omega_P$ such that $\gamma=\sum x_i(\omega)$.

{\bf Remark.}  The function $\min (x, 1)$ is chosen only for concreteness: any other positive
bounded function on $(0, +\infty)$ coinciding with $x$ on some interval $(0, \varepsilon)$
and bounded away from zero on its complement, could have been chosen instead.

Consider the set
\begin{equation*}
{\mathfrak s}\Omega_P=\biggl\{\eta\in \Mfin\bigl((0,+\infty)\bigr): \eta=\sum\limits_{i=1}^{\infty}\min (x_i, 1)\delta_{x_i}
\text{ for some }x_i>0\biggr\}.
\end{equation*}
The set ${\mathfrak s}\Omega_P$ is clearly closed in $\Mfin\bigl((0,+\infty)\bigr)$.

Any  measure $\eta$ from the set ${\mathfrak s}\Omega_P$ admits a unique representation
$\eta={\mathfrak s}\omega$ for a unique $\omega\in\Omega_P^0$.

Consequently, to any finite Borel  measure $\Prob\in \Mfin
(\Mfin((0,+\infty)))$ supported on the set ${\mathfrak s}\Omega_P$
there corresponds a unique measure
${\mathfrak p}\Prob$ on $\Omega_P$ such that ${\mathfrak s}_*{\mathfrak p}\Prob=\Prob$ and
${\mathfrak p}\Prob(\Omega_P\setminus \Omega_P^0)=0$.

\subsection{Weak convergence in $\Mfin\Mfin\left((0,+\infty)\right)$ and in $\Mfin(\Omega_P)$}

The connection of the weak convergence in the space of finite measures on the space of finite measures on the half-line
to weak convergence on the space of measures on the Pickrell set is now given
by the following
\begin{proposition}\label{weakconfpick}
Let $\nu_n, \nu\in\Mfin\Mfin\left((0,+\infty)\right)$ be supported on the set ${\mathfrak s}\Omega_P$
and assume that $\nu_n\to \nu$ weakly in $\Mfin\Mfin\left((0,+\infty)\right)$ as $n\to\infty$. Then
${\mathfrak p} \nu_n \to {\mathfrak p} \nu$ weakly in $\Mfin(\Omega_P)$ as $n\to\infty$.
\end{proposition}

The map  ${\mathfrak s}$ is, of course, not continuous, since the function $$\displaystyle\omega\to\sum\limits_{i=1}^{\infty}\min (x_i(\omega), 1)$$
is not continuous on the Pickrell set.

Nonetheless, we have the following relation between tightness of measures on $\Omega_P$ and on
$\Mfin\bigl((0,+\infty)\bigr)$.
\begin{lemma}\label{tight-pickrell}
Let $\Prob_{\alpha}\in \Mfin
(\Mfin((0,+\infty)))$ be a tight family of measures.
Then the family  ${\mathfrak p}\Prob_{\alpha}$ is also tight.
\end{lemma}

Proof. Take $R>0$ and consider the subset
$$\Omega_P(R)=\biggl\{\omega\in\Omega_P\colon \gamma(\omega)\leq R,
\sum\limits_{i=1}^{\infty}\min (x_i(\omega), 1)
\leq R\biggr\}.$$

The subset $\Omega_P(R)$ is compact in $\Omega_P$, and any compact subset of ${\mathfrak s}\Omega_P$ is
in fact a subset of ${\mathfrak s}\Omega_P(R)$ for a sufficiently large $R$.
Consequently, for any $\varepsilon>0$ one can find a sufficiently large $R$ in such a way that
$$
\Prob_{\alpha}({\mathfrak s}(\Omega_P(R)))>1-\varepsilon\text{ for all } \alpha.
$$
Since all measures $\Prob_{\alpha}$ are supported on $\Omega_P^0$, it follows that
$$
{\mathfrak p}\Prob_{\alpha}(\Omega_P(R))>1-\varepsilon\text{ for all }\alpha,
$$
and the desired tightness is established.

\begin{corollary}
Let
$$
\Prob_n\in \Mfin\left(\Mfin\bigl((0,+\infty)\bigr)\right),
n\in {\mathbb N}, \quad \Prob\in \Mfin\left(\Mfin\bigl((0,+\infty)\bigr)\right)
$$
be finite Borel measures.
Assume \begin{enumerate}
\item the measures $\Prob_n$ are supported on the set ${\mathfrak s}\Omega_P$
for all $n\in\mathbb{N}$;
\item $\Prob_n\to \Prob$ converge weakly in $\Mfin\left(\Mfin\bigl((0,+\infty)\bigr)\right)$  as $n\to\infty$.
\end{enumerate}
Then the measure $\Prob$ is also supported on   the set ${\mathfrak s}\Omega_P$
and  ${\mathfrak p}\Prob_n\to{\mathfrak p}\Prob$ weakly in $\Mfin\left(\Omega_P\right)$ as $n\to\infty$.
\end{corollary}

{Proof}. The measure $\Prob$ is of course supported on   the set ${\mathfrak s}\Omega_P$, since the set ${\mathfrak s}\Omega_P$
 is closed. The desired  weak convergence in $\Mfin\left(\Omega_P\right)$
is now established in three steps.

{\it The First Step: The Family ${\mathfrak p}\Prob_n$ is Tight.}

The family ${\mathfrak p}\Prob_n$ is tight by Lemma \ref{tight-pickrell} and therefore admits
a weak accumulation point $\Prob'\in \Mfin\left(\Omega_P\right)$.

{\it The Second Step: Finite-Dimensional Distributions Converge.}

Let $l\in\mathbb{N}$, let $\varphi_l\colon(0,+\infty)\to\mathbb{R}$ be continuous compactly supported functions, set $\varphi_l(x)=\min (x, 1)\psi_l(x)$, take $t_1,\dots,t_l\in\mathbb{R}$ and observe that, by definition, for any $\omega\in\Omega_P$ we have
\begin{equation*}
\exp\left(i\sum\limits_{k=1}^lt_k\left(\sum\limits_{r=1}^{\infty}\varphi_k\left(x_r(\omega)\right)\right)\right)= \exp\left(i\sum\limits_{k=1}^lt_k\Int_{\psi_k}\left({\mathfrak s}\omega\right)\right)
\end{equation*}
and consequently
\begin{multline*}
\int\limits_{\Omega_P}\exp\left(i\sum\limits_{k=1}^lt_k\left(\sum\limits_{r=1}^{\infty}\varphi_k\left(x_r(\omega) \right)\right)\right)\,d\Prob'(\omega)=\\
=\int\limits_{\Mfin((0,+\infty))} \exp\left(i\sum\limits_{k=1}^lt_k\Int_{\psi_k}(\eta)\right)\,d({\mathfrak s}_*\Prob'(\eta)).
\end{multline*}

We now write
\begin{multline*}
\int\limits_{\Omega_P}\exp\left(i\sum\limits_{k=1}^lt_k\left(\sum\limits_{r=1}^{\infty}\varphi_k\left(x_r(\omega) \right)\right)\right)\,d\Prob'(\omega)=\\
=\lim\limits_{n\to\infty}\int\limits_{\Omega_P}\exp\left(i\sum\limits_{k=1}^lt_k\left(\sum\limits_{r=1}^{\infty} \varphi_k\left(x_r(\omega)\right)\right)\right)\,d\Prob_n(\omega)\,.
\end{multline*}
On the other hand, since $\Prob_n\to\Prob$ weakly in $\Mfin\left(\Mfin\bigl((0,+\infty)\bigr)\right)$, we have
\begin{multline*}
\lim\limits_{n\to\infty}\int\limits_{\Mfin\bigl((0,+\infty)\bigr)} \exp\left(i\sum\limits_{k=1}^lt_k\Int_{\psi_k}(\eta)\right)\,d({\mathfrak s}_*\Prob_n)=\\
=\int\limits_{\Mfin((0,+\infty))} \exp\left(i\sum\limits_{k=1}^lt_k\Int_{\psi_k}(\eta)\right)\,d\Prob.
\end{multline*}

It follows that
\begin{multline*}
\int\limits_{\Omega_P}\exp\left(i\sum\limits_{k=1}^lt_k\left(\sum\limits_{r=1}^{\infty}\varphi_k\left(x_r(\omega) \right)\right)\right)\,d\Prob'(\omega)=\\
=\int\limits_{\Mfin\bigl((0,+\infty)\bigr)} \exp\left(i\sum\limits_{k=1}^lt_k\Int_{\psi_k}(\eta)\right)\,d\Prob.
\end{multline*}

Since integrals of functions of the form $\exp\left(i\sum\limits_{k=1}^lt_k\Int_{\psi_k}(\eta)\right)$ determine
a finite Borel measure on $\Mfin\bigl((0,+\infty)\bigr)$ uniquely, we have
$${\mathfrak s}_*\Prob'=\Prob\,.$$

{\it The Third Step: The Limit Measure is Supported on $\Omega_P^0$.}

To see that $\Prob'(\Omega_P\backslash\Omega_P^0)=0$, set
$$
\gamma'(\omega)=\gamma(\omega)+\sum\limits_{k:x_k(\omega)\geq 1} (1-x_k(\omega)).
$$
Since the sum in the right-hand side is finite, the function $\gamma'$ is continuous on $\Omega_P$, and we have
$$
\int\limits_{\Omega_P}\exp({-\gamma'(\omega)})\,d\Prob'(\omega)= \lim\limits_{n\to\infty}\int\limits_{\Omega_P}\exp({-\gamma'(\omega)})\,d\mathfrak{p}\Prob_n(\omega).
$$
We also have
\begin{multline*}
\int\limits_{\Omega_P}\exp\biggl(-\sum\limits_{i=1}^{\infty}\min (1, x_i(\omega))\biggr)\,d\Prob'(\omega)={}\\  \lim\limits_{n\to\infty}\int\limits_{\Omega_P}\exp\biggl(-\sum\limits_{i=1}^{\infty}\min (1, x_i(\omega))\biggr)\,d\mathfrak{p}\Prob_n(\omega).
\end{multline*}
Since  $\mathfrak{p}\Prob_n(\Omega_P\backslash\Omega_P^0)=0$ for all $n$, we have
$$
\int\limits_{\Omega_P}\exp({-\gamma'(\omega)})\,d\mathfrak{p}\Prob_n(\omega)= \int\limits_{\Omega_P}\exp\biggl(-\sum\limits_{i=1}^{\infty}\min (1, x_i(\omega))\biggr)\,d\mathfrak{p}\Prob_n(\omega),
$$
for all $n$. It follows that
$$
\int\limits_{\Omega_P}\exp(-\gamma'(\omega))\,d\Prob'(\omega)= \int\limits_{\Omega_P}\exp\biggl(-\sum\limits_{i=1}^{\infty}\min (1, x_i(\omega))\biggr)\,d\Prob'(\omega),
$$
whence the equality $\gamma'(\omega)=\sum\limits_{i=1}^{\infty}\min (1, x_i(\omega))$
and, consequently, also the equality $\gamma(\omega)=\sum\limits_{i=1}^{\infty}x_i(\omega)$ holds $\Prob'$-almost surely, and $\Prob'(\Omega_P\backslash\Omega_P^0)=0$.
We thus have  $\Prob'={\mathfrak p}\Prob$. The proof is complete.

\section{Proof of Lemma \ref{main-lemma} and Completion of the Proof of Theorem \ref{mainthm}}

\subsection{Reduction of Lemma \ref{main-lemma} to Lemma \ref{fast-decay}}

Recall that we have introduced a sequence of mappings
$$\rrad^{(n)} : \Mat(n, \mathbb{C}) \rightarrow \Omega_P^0, \quad n \in \mathbb{N}$$
that to $z \in \Mat(n, \mathbb{C})$
assigns the point
$$\rrad^{(n)}(z)=\left( \frac{\tr z^* z}{n^2}, \frac{\lambda_1(z)}{n^2}, \dots ,
\frac{\lambda_n (z)}{n^2}, 0, \dots, 0, \dots \right),$$
where
$\lambda_1 (z)\geq \dots \geq \lambda_n (z) \geq
0$ are the eigenvalues of the matrix $z^* z$, counted with
multiplicities and arranged in non-increasing order.
By definition, we have
$$
\gamma(\rrad^{(n)}(z))=\frac{\tr z^* z}{n^2}.
$$
Following Vershik \cite{Vershik},
we now introduce on  $\Mat(\mathbb{N}, \mathbb{C})$ a sequence of averaging operators over the compact groups $U(n)\times U(n)$.
\begin{equation*}
\left( \mathcal A_n f\right) (z)=\displaystyle \int\limits_{U(n)\times U(n)} f(u_1zu_2^{-1})du_1 du_2,
\end{equation*}
where $du$ stands for the normalized Haar measure on the group $U(n)$.
For any $U(\infty)\times U(\infty)$-invariant probability measure on  $\Mat(\mathbb{N}, \mathbb{C})$,
the operator $\mathcal A_n$ is the operator of conditional  expectation with respect to the sigma-algebra of $U(n)\times U(n)$-invariant sets.

By definition,  the function $\left( \mathcal A_n f\right) (z)$ only depends on $\rrad^{(n)}(z)$.

\begin{lemma}\label{fast-decay}
Let $m\in {\mathbb N}$.
There exists a positive Schwartz function $\varphi$ on $\Mat(m, \mathbb{C})$ as well as
a positive continuous function
$f$ on $\Omega_P$ such that for any $z\in \Mat(\mathbb{N}, \mathbb{C})$ and any $n\geq m$
we have
\begin{equation*}
f(\rrad^{(n)}(z))\leq \left( \mathcal A_n \varphi\right) (z).
\end{equation*}
\end{lemma}

{\bf Remark.} The function $\varphi$, initially defined on $\Mat(m, \mathbb{C})$, is here extended to
$\Mat(\mathbb{N}, \mathbb{C})$ in the obvious way: the value of $\varphi$ at a matrix $z$ is set to be its value on its $m\times m$-corner.

We postpone the proof of the Lemma to the next subsection and proceed with the the proof of
Lemma \ref{main-lemma}.

Refining the definition of the class $\mathfrak{F}$ in the introduction to the first part of the paper, take $m\in {\mathbb N}$ and let $\mathfrak{F}(m)$ be the family of all Borel sigma-finite
$U(\infty)\times U(\infty)$-invariant measures $\nu$
on $\Mat(\mathbb{N}, \mathbb{C})$ such that for any $R>0$  we have
$$
\nu\left(\{z: \max\limits_{i,j\leq m} |z_{ij}|<R\}\right)<+\infty.
$$
Equivalently, the measure of a set of matrices, whose $m\times m$-corners are required to lie
in a compact set, must be finite; in particular, the projections $(\pi^{\infty}_n)_*\nu$
are well-defined   for all $n\ge m$.
For example, if $s+m>0$, then the Pickrell measure $\mu^{(s)}$ belongs to $\mathfrak{F}(m)$.
Recall furthermore that, by the results of \cite{Buf-ergdec}, \cite{Buf-inferg}  any measure $\nu\in \mathfrak{F}(m)$ admits a unique ergodic decomposition into {\it finite} ergodic components: in other words, for any such $\nu$ there exists
a unique Borel sigma-finite measure ${\overline \nu}$ on $\Omega_P$ such that we have
\begin{equation*}
\nu=\int\limits_{\Omega_P}\eta_{\omega} d{\overline \nu}(\omega).
\end{equation*}

Since the orbit of the unitary group is of course a compact set, the measures $(\rrad^{(n)})_*\nu$ are well-defined for $n>m$ and may be thought of as finite-dimensional approximations of the decomposing measure ${\overline \nu}$. Indeed, recall from the introduction to the first part of the paper that, if $\nu$ is finite, then the measure
${\overline \nu}$ is the weak limit of the measures $(\rrad^{(n)})_*\nu$ as $n\to\infty$.
The following proposition is a stronger and a more precise version of Lemma \ref{main-lemma} from the introduction.
\begin{proposition}
Let $m\in {\mathbb N}$, let  $\nu\in \mathfrak{F}(m)$,
let $\varphi$ and $f$ be given by Lemma \ref{fast-decay}, and assume
$$
\varphi\in L_1(\Mat(\mathbb{N}, \mathbb{C}), \nu).
$$
Then
\begin{enumerate}
\item $f\in L_1(\Omega_P,(\rrad^{(n)})_* {\nu})$ for all $n>m$;
\item $f\in L_1(\Omega_P, {\overline \nu})$;
\item $f(\rrad^{(n)})_* {\nu}\to f{\overline \nu}$ weakly in $\Mfin(\Omega_P)$.
\end{enumerate}
\end{proposition}

Proof.
{\it First Step: The Martingale Convergence Theorem and the Ergodic Decomposition.}

We start by formulating a pointwise version of the equality (\ref{rinftymu})  from the Introduction:
for any $z\in \Matreg$ and any bounded continuous function $\varphi$ on
$\Mat(\mathbb{N};\mathbb{C})$
we have
\begin{equation}\label{anrinfty}
\lim\limits_{n\to\infty} {\mathcal A}_n\varphi(z)=\int\limits_{\Mat(\mathbb{N};\mathbb{C})} fd\eta_{\rrad^{\infty}(z)}
\end{equation}
(here, as always, given $\omega\in\Omega_P$,   the symbol $\eta_{\omega}$ stands for the ergodic probability measure corresponding to $\omega$.)
Indeed, (\ref{anrinfty}) immediately follows from the definition of regular matrices, the  Olshanski-Vershik characterization of the convergence of orbital measures \cite{OlshVershik} and   the Reverse Martingale Convergence Theorem.

{\it The Second Step.}

Now  let $\varphi$ and $f$ be given by Lemma \ref{fast-decay}, and assume
$$
\varphi\in L_1(\Mat(\mathbb{N}, \mathbb{C}), \nu).
$$
\begin{lemma}For any $\varepsilon >0$ there exists a
$U(\infty) \times U(\infty)$-invariant set
 $Y_{\varepsilon}
\subset \Mat(\mathbb{N}, \mathbb{C})$ such that
\begin{enumerate}
\item $\nu(Y_\varepsilon) < + \infty ;$
\item for all  $n>m$ we have
$$\int \limits _{\Mat(\mathbb{N}, \mathbb{C})\setminus Y_{\varepsilon}} f(\rrad^{(n)}(z))
d\nu(z) < \varepsilon.$$
\end{enumerate}
\end{lemma}
Proof.
Since $\varphi\in L_1(\Mat(\mathbb{N}, \mathbb{C}), \nu)$, we have
$$
\int\limits_{\Omega_P}\Biggl(\int\limits_{\Mat(\mathbb{N}, \mathbb{C})}\varphi d\eta_{\omega}\Biggr) d{\overline \nu}(\omega)<+\infty.
$$

Choose a Borel subset  ${\tilde Y}^{\varepsilon}\subset \Omega_P$  in such a way that
${\overline \nu}({\tilde Y}^{\varepsilon})<+\infty$ and
$$
\int\limits_{{\tilde Y}^{\varepsilon}}\Biggl(\int\limits_{\Mat(\mathbb{N}, \mathbb{C})}\varphi d\eta_{\omega}\Biggr) d{\overline \nu}(\omega)<\varepsilon.
$$
The pre-image of the set  ${\tilde Y}^{\varepsilon}$ under the map $\rrad^{(\infty)}$ or, more precisely, the set
$$
Y_{\varepsilon}=\{z\in \Matreg: \rrad^{(\infty)}(z)\in {\tilde Y}^{\varepsilon} \}
$$
is by definition $U(\infty) \times U(\infty)$-invariant and has all the desired properties.

{\it The Third Step.}

Let $\psi : \Omega_P \rightarrow \mathbb{R}$ be
continuous and bounded. Take $\varepsilon >0$ and the corresponding set $Y_{\varepsilon}$.

For any $z\in\Matreg$ we have
$$
\lim \limits _{n \to \infty}
\psi(\rrad^{(n)}(z))
\cdot f(\rrad^{(n)}(z))=\psi(\rrad^{(\infty)}(z)) \cdot f(\rrad^{(\infty)}(z)) .
$$
Since $\nu(Y_{\varepsilon})<\infty$, the bounded convergence theorem gives
\begin{multline*}
\lim \limits _{n \to \infty} \int \limits _{Y_{\varepsilon}}
\psi(\rrad^{(n)}(z))
\cdot f(\rrad^{(n)}(z)) d\nu(z)=\\
=\int \limits _{Y_{\varepsilon}} \psi(\rrad^{(\infty)}(z)) \cdot f(\rrad^{(\infty)}(z)) d\nu(z).
\end{multline*}
By definition of $Y_{\varepsilon}$
for all $n \in \mathbb{N}$, $n>m$,  we have
$$\left| \int \limits _{\Mat(\mathbb{N}, \mathbb{C})
 \setminus Y_{\varepsilon}} \psi(\rrad^{(n)}(z)) \cdot f(\rrad^{(n)}(z)) d\nu(z) \right|
 \leq \varepsilon \sup\limits_{\Omega_P} |\psi|.$$
It follows that
\begin{multline*}
\lim\limits_{n\to\infty} \int\limits_{\Mat(\mathbb{N}, \mathbb{C})}
\psi(\rrad^{(n)}(z))
\cdot f(\rrad^{(n)}(z)) d\nu(z)=\\
=\int\limits_{\Mat(\mathbb{N}, \mathbb{C})}
\psi(\rrad^{(\infty)}(z)) \cdot f(\rrad^{(\infty)}(z)) d\nu(z),
\end{multline*}
which, in turn, implies that
$$
\lim\limits_{n \to \infty} \int\limits_{\Omega_P}
\psi f d(\rrad^{(n)})_*(\nu)=
\int\limits_{\Omega_P} \psi f d\overline{\nu},
$$
 that the weak convergence is established, and that the Lemma is proved completely.

\subsection{Proof of Lemma \ref{fast-decay}}

Introduce an inner product $\langle\,\cdot\,,\,\cdot\,\rangle$ on $\Mat(m, \mathbb{C})$ by the formula
$\langle z_1, z_2\rangle=\Re \tr(z_1^*z_2)$.
This inner product is naturally extended to a pairing between the projective limit
$\Mat({\mathbb N}, \mathbb{C})$ and
the inductive limit
$$
\Mat_0=\bigcup\limits_{m=1}^{\infty} \Mat(m, \mathbb{C}).
$$
For a matrix $\zeta\in \Mat_0$ set
$$
\Xi_{\zeta}(z)=\exp (i\langle \zeta, z \rangle), \ z\in \Mat({\mathbb N}, \mathbb{C}).
$$
% as a function of
%$z\in \Mat({\mathbb N}, \mathbb{C})$.

We start with the following simple estimate on the behaviour  of the
Fourier transform of orbital measures.

\begin{lemma}\label{est-fourier}
Let $m\in {\mathbb N}$.
For any $\varepsilon>0$ there exists $\delta>0$ such that for any $n>m$ and
$\zeta\in \Mat(m, \mathbb{C})$, $z\in  \Mat({\mathbb N}, \mathbb{C})$ satisfying
$$
\tr(\zeta^*\zeta) \tr \left((\pi_n^{\infty}(z))^*(\pi_n^{\infty}(z)\right)<\delta n^2
 $$
we have
$$
|1-{\mathcal A}_n\Xi_{\zeta}(z)|<\varepsilon.
$$
\end{lemma}
Proof.
This is a simple corollary of the power series representation of the Harish-Chandra--Itzykson--Zuber orbital integral, see e.g. \cite{Faraut}, \cite{Faraut-angl}, \cite{Rabaoui2}. Indeed, let $\sigma_1, \dots, \sigma_m$ be the eigenvalues of $\zeta^*\zeta$, and let $x_1^{(n)}, \dots, x_{n}^{(n)}$ be the eigenvalues of $\pi_n^{\infty}(z)$.

The standard power series representation, see e.g.  \cite{Faraut}, \cite{Faraut-angl},\cite{Rabaoui2}, for the Harish-Chandra--Itzykson--Zuber orbital integral gives, for any $n \in \mathbb{N}$, a representation
$$
A_n \Xi_\zeta(z) = 1 + \sum\limits_{\lambda \in \mathbb{Y}_+} a(\lambda, n) s_{\lambda}\left(\sigma_1, \dots, \sigma_m\right) \cdot s_{\lambda}\left(\frac{x_1^{(n)}}{n^2}, \dots, \frac{x_n^{(n)} }{n^2}\right),
$$
where the summation takes place over the set $\mathbb{Y}_+$ of all non-empty Young diagrams $\lambda$, $s_{\lambda}$ stands for the Schur polynomial corresponding to the diagram $\lambda$, and the coefficients $a(\lambda, n)$ satisfy
$$
\sup\limits_{\lambda \in \mathbb{Y}_+} |a(\lambda, n) | \leq 1
$$
The lemma follows immediately.

\begin{corollary}\label{schw-R}
For any $m\in {\mathbb N}$, $\varepsilon>0$, $R>0$,  there exists a positive Schwartz function
$\psi:\Mat(m, \mathbb{C})\to (0,1]$ such that for all $n>m$ we  have
\begin{equation}\label{schw-ineq}
{\mathcal A}_n\psi(\pi_m^{\infty}(z))\geq 1-\varepsilon
\end{equation}
for all $z$ satisfying
$$
\tr \left((\pi_n^{\infty}(z))^*(\pi_n^{\infty}(z)\right)<Rn^2.
$$
\end{corollary}

Proof. Let $\psi$ be a Schwartz function  taking values in $(0,1]$. Assume additionally that
$\psi(0)=1$ and  that the Fourier transform  of $\psi$ is supported in
the ball of radius $\varepsilon_0$ around the origin. A Schwartz function satisfying all these requirements is constructed without difficulty.
By Lemma \ref{est-fourier}, if $\varepsilon_0$ is small enough as a function of $m, \varepsilon, R$, then
the  inequality (\ref{schw-ineq}) holds for all $n>m$. Corollary \ref{schw-R} is proved completely.

We now conclude the proof of Lemma \ref{fast-decay}.

Take a sequence $R_n\to\infty$, and let $\psi_n$ be the corresponding sequence of Schwartz functions
given by Corollary \ref{schw-R}.
Take positive numbers $t_n$ decaying fast enough so that the function
$$
\varphi=\sum_{n=1}^{\infty} t_n\psi_n
$$
is Schwartz.

Let ${\tilde f}$ be a positive continuous function on $(0, +\infty)$ such that for any $n$,  if $t\leq R_n$,
then ${\tilde f}(t)<t_n/2$.
For $\omega\in\Omega_P$, $\omega=(\gamma, x)$, set
$$
f(\omega)={\tilde f}(\gamma(\omega)).
$$

The function $f$ is by definition positive and continuous.   By Corollary \ref{schw-R}, the functions $\varphi$ and $f$ satisfy all requirements of Lemma \ref{fast-decay}, which, therefore, is proved completely.

\subsection {Completion of the proof of Theorem \ref{mainthm}}
\begin{lemma}\label{fgpb}
Let $E$ be a locally compact complete metric space.
Let $\mathbb B_n, \mathbb B$ be sigma-finite measures on $E$, let $\Prob$ be a probability measure on $E$,
and let $f, g$ be positive bounded continuous functions on $E$.
Assume that for all $n\in {\mathbb N}$ we have
$$
g\in L_1(E, \mathbb B_n)
$$
and that, as $n\to\infty$, we have
\begin{enumerate}
\item $\displaystyle f\mathbb B_n\to f\mathbb B$ weakly in $\Mfin(E)$;
\item $\displaystyle \frac{g\mathbb B_n}{\displaystyle \int\limits_E gd\mathbb B_n} \to \Prob$ weakly in $\Mfin(E)$.
\end{enumerate}
Then
$$
g\in L_1(E, \mathbb B)
$$
and
$$
\Prob=\frac{g\mathbb B}{\displaystyle \int\limits_E gd\mathbb B}.
$$
\end{lemma}
Proof. Let $\varphi$ be a nonnegative bounded continuous function on $E$. On the one hand, as $n\to\infty$, we have
$$
\displaystyle \int\limits_E \varphi f g d{\mathbb B}_n\to \displaystyle \int\limits_E \varphi f g d{\mathbb B},
$$
and, on the other  hand, we have
\begin{equation}\label{phifg}
\frac{\displaystyle \int\limits_E \varphi f g d{\mathbb B}_n}{\displaystyle \int\limits_E gd\mathbb B_n}\to \displaystyle  \int\limits_E \varphi f  d\Prob.
\end{equation}
Choosing $\varphi=1$, we obtain that
$$
\lim\limits_{n\to\infty}\int\limits_E gd\mathbb B_n=\frac{\int\limits_E  f g d{\mathbb B}}{\int\limits_E  f  d\Prob}>0;
$$
the sequence $\displaystyle \int\limits_E gd\mathbb B_n$ is thus bounded away both from zero and infinity. Furthermore, for arbitrary bounded continuous positive $\varphi$ we have
\begin{equation*}
\lim\limits_{n\to\infty}\displaystyle \int\limits_E gd\mathbb B_n=\displaystyle \frac{\displaystyle \int\limits_E \varphi f g d{\mathbb B}}{\displaystyle \int\limits_E \varphi f  d\Prob}.
\end{equation*}
Now take $R>0$ and  $\varphi(x)=\min (1/f(x), R)$. Letting $R$ tend to $\infty$, we obtain
\begin{equation}\label{gdb}
\lim\limits_{n\to\infty}\displaystyle \int\limits_E gd\mathbb B_n={\displaystyle \int\limits_E  g d{\mathbb B}}.
\end{equation}
Substituting (\ref{gdb}) back into (\ref{phifg}), we arrive at the equality
$$
\displaystyle \int\limits_E \varphi f  d\Prob=\frac{\displaystyle \int\limits_E \varphi f g d{\mathbb B}}{\displaystyle \int\limits_E  g d{\mathbb B}}.
$$
Note that here, as in (\ref{phifg}), the function $\varphi$ may be an arbitrary nonnegative continuous function on $E$. In particular,  taking a compactly supported function
$\psi$ on $E$  and setting $\varphi=\psi/f$, we obtain
$$
\displaystyle \int\limits_E \psi d\Prob=\displaystyle \frac{\displaystyle \int\limits_E \psi g d{\mathbb B}}{\displaystyle \int\limits_E  g d{\mathbb B}}.
$$
Since this equality is true for any compactly supported fnction $\psi$ on $E$, we conclude that
$$
\Prob=\frac{g\mathbb B}{\displaystyle \int\limits_E gd\mathbb B},
$$
and the Lemma is proved completely.

Combining Lemma \ref{fgpb} with Lemma \ref{main-lemma} and Proposition \ref{weak-pick}, we conclude the proof of Theorem \ref{mainthm}.

Theorem \ref{mainthm} is proved completely.

\subsection {Proof of Proposition \ref{singul}}
Using Kakutani's theorem, we  now conclude the proof of Proposition \ref{singul}.
Take $n$ large enough so that $n+s>1$, $n+s'>1$ and compute the Hellinger integral
\begin{multline*}
\mathrm{Hel}\,(n,s,s')=\mathbb{E}\left(\sqrt{\left(P^{(n,n-1,s)}\!\!\times\!P^{(n,n,s)}\right) \!\cdot\!\left(P^{(n,n-1,s')}\!\!\times\!P^{(n,n,s')}\right)}\right)= \\
=\sqrt{\frac{\Gamma(2n-1+s)}{\Gamma(n-1)\Gamma(n+s)} \cdot \frac{\Gamma(2n-1+s')}{\Gamma(n-1)\Gamma(n+s')}\cdot \frac{\Gamma(2n+s)}{\Gamma(n)\Gamma(n+s)}  \cdot
\frac{\Gamma(2n+s')}{\Gamma(n)\Gamma(n+s')} }
\times \\
\times\int\limits_0^{\infty}r^{n-1}(1+r)^{-2n-1-\frac{s+s'}{2}}\!dr\!\cdot\! \int\limits_0^{\infty}r^{n-1}(1+r)^{-2n-\frac{s+s'}{2}}\!dr= \\
=\frac{\sqrt{\Gamma(2n-1+s)\!\cdot\!\Gamma(2n-1+s')}}{\Gamma(2n-1+\frac{s+s'}{2})}\cdot \frac{\sqrt{\Gamma(2n+s)\!\cdot\!\Gamma(2n+s')}}{\Gamma(2n+\frac{s+s'}{2})}\cdot \frac{\left(\Gamma(n+\frac{s+s'}{2})\right)^2}{\Gamma(n+s)\!\cdot\!\Gamma(n+s')}\,.\notag
\end{multline*}

We now recall a classical asymptotics: as $t\to\infty$, we have
$$\frac{\Gamma(t+a_1)\!\cdot\!\Gamma(t+a_2)}{\left(\Gamma(t+\frac{a_1+a_2}{2})\right)^2}= 1+\frac{(a_1+a_2)^2}{4t}+O\left(\frac{1}{t^2}\right).$$
It follows that
$$\mathrm{Hel}(n,s,s')=1-\frac{(s+s')^2}{8\,n}+O\left(\frac{1}{n^2}\right)\,,$$
whence, by the Kakutani's theorem combined with (\ref{inf-prod-s}), the
Pickrell measures $\mu^{(s)}$ and $\mu^{(s')}$, finite or infinite,  are mutually
singular if $s\neq s'$.

\subsection{Proof of Proposition \ref{bs-sing}}

In view of Proposition \ref{mainthm-fin} and Theorem \ref{mainthm}, it suffices to prove the singularity
of the ergodic decomposition measures ${\overline \mu}^{(s_1)}$, ${\overline \mu}^{(s_2)}$. Since, by Proposition \ref{singul},  the measures
 ${ \mu}^{(s_1)}$, ${ \mu}^{(s_2)}$ are mutually singular, there exists a set $D\subset \Mat(\mathbb{N}, \mathbb{C})$
such that
$${\mu}^{(s_1)}(D)=0,\quad { \mu}^{(s_2)}(\Mat(\mathbb{N}, \mathbb{C})\setminus D)=0.$$
Introduce the set
$$
{\tilde D}=\{z\in  \Mat(\mathbb{N}, \mathbb{C}): \lim\limits_{n\to\infty} {\mathcal A}_n\chi_D(z)=1\}.
$$
By definition, the set ${\tilde D}$ is $U(\infty)\times U(\infty)$-invariant, and
we have  $${ \mu}^{(s_1)}({\tilde D})=0,\quad
{\mu}^{(s_2)}( \Mat(\mathbb{N}, \mathbb{C})\setminus {\tilde D})=0.$$
Introduce now the set ${\overline D}\subset \Omega_P$ by the formula
$$
{\overline D}=\{\omega\in \Omega_P: \eta_{\omega}({\tilde D})=1\}.
$$
We clearly have
 $${\overline \mu}^{(s_1)}({\overline D})=0, \quad {\overline \mu}^{(s_2)}( \Omega_P\setminus {\overline D})=0.$$

Proposition \ref{bs-sing} is proved completely.

{\bf Acknowledgements.}
This project has received funding from the European Research Council (ERC) under the European Union's Horizon 2020 research and innovation programme (grant agreement No 647133 (ICHAOS))
and has also been funded by Grant MD 5991.2016.1 of the President of the Russian Federation as well as by 
the Russian Academic Excellence Project `5-100'.

\end{document}